\documentclass[a4paper,twoside,11pt]{article}
\usepackage{amsfonts}
\usepackage{amssymb}\usepackage{amsmath}
\usepackage{graphicx}
\def\be{\begin{equation}}
\def\ee{\end{equation}}
\def\ba{\begin{array}}
\def\ea{\end{array}}
\def\bc{\begin{center}}
\def\ec{\end{center}}

\def\ZZ{\rm {{\rm Z}\kern-.48em{\rm Z}}}
\def\RR{\rm \hbox{I\kern-.2em\hbox{R}}}
\def\CC{\rm \hbox{C\kern -.5em{\raise .32ex
\hbox{$\scriptscriptstyle |$}}\kern - .22em{\raise .6ex
\hbox{$\scriptscriptstyle |$}}\kern
.4em}}
\def\CC{ {\mathbf C} }
\def\RR{ {\mathbf R} }
\def\ZZ{ {\mathbf Z} }

\newtheorem{Tm}{Theorem}[section]
\newtheorem{Df}[Tm]{Definition}%[section]

\newtheorem{Lm}[Tm]{Lemma}%[section]
 \numberwithin{equation}{section}
\newcommand{\double}{\baselineskip 1.24 \baselineskip}

\title{ Random sampling stability in weighted reproducing kernel subspaces of  $L_\nu^p(\mathbb{R}^d)$
}
\author{{\small  Yingchun Jiang,\ \ Yajing Zhang, \ \ Wan Li}\\
{\small
 School of Mathematics and Computational Science,}\\{\small
Guilin University of Electronic Technology, Guilin,  P. R. China}
}

\begin{document}
\date{}
\maketitle \double
\textbf{Abstract:}\ In this paper, we mainly study the random sampling stability  for signals in  a weighted reproducing kernel subspace
 of $L_\nu^p(\mathbb{R}^d)$ without the additional requirement that the kernel function has symmetry. The sampling set is independently and randomly drawn from a general probability distribution over  $\mathbb{R}^d$. Based on the frame characterization of weighted reproducing kernel subspaces,  we first approximate the weighted reproducing kernel space by a finite dimensional subspace  on any bounded domains. Then, we prove that the random sampling stability holds with high probability for all signals in weighted reproducing kernel subspaces  whose energy concentrate on a cube when the sampling size is large enough.


 {\bf MR(2000) Subject Classification:}  94A20, 46E30.

\section{Introduction }

\ \ \ \ Random sampling plays an important role in many fields, such as image processing \cite{C2}, compressed sensing \cite{E1} and learning theory \cite{S2}. Random sampling has been generally studied for multivariate trigonometric polynomials \cite{B4}, bandlimited signals \cite{B1,B2}, signals
that satisfy some locality properties in short-time Fourier transform \cite{V1}, signals with bounded derivatives \cite{Y2}, signals in a shift-invariant space \cite{F1,L2,Y1,Y3}, signals with finite rate of innovation \cite{L3} and signals in  reproducing kernel subspaces of $L^p(\mathbb{R}^d)$ \cite{L1,P1}.

Stability and reconstruction algorithm are two fundamental problems in sampling theory. In \cite{P1}, sampling stability was established with high probability for signals in  energy concentrated subspaces of reproducing kernel spaces. Because such subspaces are nonlinear and almost all reconstruction algorithms  were only given for  functions in a finite dimensional subspace \cite{P1,Y1,Y3}, an iterative algorithm which provides approximation to signals with energy concentrated on a cube was firstly constructed in \cite{L1}. Note that random samples  in \cite{P1} were taken from a uniform distribution on a bounded domain and the kernel function $K$ was assumed to
satisfy a very strong symmetric condition
\begin{equation}\label{con1}
%|K(x,y)|\leq \frac{C}{(1+\|x-y\|_1)^\alpha},\ \alpha>d(1-1/p)+d+1\ \  {\rm{and}} \ \
 K(x,y)=K(y,x).
\end{equation}
In this paper, we will restudy the  random sampling stability for signals in a weighted reproducing kernel subspace of $L_\nu ^p(\mathbb{R}^d)$ without the additional condition %$K(x,y)=K(y,x)$ in
\eqref{con1}. Moreover, the random samples are drawn over $\mathbb{R}^d$ from a general probability distribution with density function $\rho$ satisfying
\begin{equation}\label{den}
0<c_\rho={\rm{ess}}\inf \limits_{x\in C_R}\rho(x)\ {\rm{and}} \ C_\rho={\rm{ess}}\sup \limits_{x\in \mathbb{R}^d}\rho(x)<\infty,
\end{equation} where $C_{R}=[-R,R]^{d}$ for $R>0$, ${\rm{ess}}\inf$ and ${\rm{ess}}\sup$ are essential infimum and supremom, respectively.
 In fact, random sampling with similar  probability distribution had been  introduced in \cite{L2,L3} for shift-invariant signals and signals with finite rate of innovation.

Suppose that $\omega$ is a weight function which is continuous, symmetric, positive and submultiplicative,
\begin{equation}
0<\omega(x+y)\leq\omega(x)\omega(y),\ x,y\in\mathbb{R}^d.
\end{equation}
Weight function $\nu$ is said to be $\omega$-moderate, that is, it is continuous, symmetric, positive and satisfies
\begin{equation}
0<\nu(x+y)\leq C_0\omega(x)\nu(y),\ x,y\in\mathbb{R}^d
\end{equation}
for some positive constant $C_0>0$. More details about weight functions can refer to \cite{G1}.

For $1\leq p\leq \infty$, $L_\nu^p(\mathbb{R}^d)$ is the Banach space of all weighted  $p$-integrable function on $\mathbb{R}^d$,
\begin{equation}
L_\nu^p(\mathbb{R}^d) = \{ f:\;{\left\| f \right\|_{L_\nu^p}} = {\left\| {\nu f} \right\|_{{L^p}}} < \infty \}.
\end{equation}

 We assume that $K(x,y)$ satisfies
 \begin{equation}\label{jz1}
|K(x,y)|\omega \left( {y - x} \right) \le \frac{{{\widetilde{C}}}}{{{{(1 + {{\left\| x \right\|}_1} + {{\left\| y \right\|}_1})}^\alpha }}},\;\alpha  > d.
\end{equation}
 %Suppose that $\omega$ is a weight function which is continuous, symmetric, positive and sub-multiplicative,
%\begin{equation}\label{j3}
%0<\omega(x+y)\leq \omega(x)\omega(y),\ x,y\in \mathbb{R}^d.
%\end{equation}
%Weight function $\nu$ is said to be $\omega$-moderate, that is, it is continuous, symmetric, positive and satisfies
%\begin{equation}\label{j1}
%0<\nu(x+y)\leq C_0\omega(x)\nu(y),\ x,y\in \mathbb{R}^d
%\end{equation}
%for some positive constant $C_0>0$. More details about weight functions can refer to \cite{G1}.
Then it is easy to verify that
\begin{align}\label{lab1}
\|K\|_{\mathcal{W}}=\max\Big\{\Big\|\sup \limits_{z\in \mathbb{R}^d}|K(z,\cdot+z)|\Big\|_{L_\omega ^1},\Big\|\sup \limits_{z\in \mathbb{R}^d}|K(\cdot+z,z)|\Big\|_{L_\omega ^1}\Big\}<\infty.
\end{align}
In fact, both the exponential kernel and the gaussian kernel satisfy the condition \eqref{jz1}. Moreover, $\|K\|_{\mathcal{W}}\leq \frac{2^d\widetilde{C}}{(\alpha-1)(\alpha-2)\cdots (\alpha-d)}$. Furthermore, we assume that
\begin{equation}\label{lab2}
\lim \limits_{\delta\rightarrow 0}\|\omega_\delta(K)\|_{\mathcal{W}}=0.
\end{equation}
Here, $\omega_\delta(K)$ is the modulus of continuity defined by
\begin{equation}
\omega_\delta(K)(x,y)=\sup \limits_{|x'|,|y'|\leq \delta}|K(x+x',y+y')-K(x,y)|.
\end{equation}
Suppose that $T$ is an idempotent $(T^2=T)$ integral operator with kernel $K$,
\begin{equation}
Tf(x)=\int_{\mathbb{R}^d}K(x,y)f(y)dy,\ f\in L_\nu^p(\mathbb{R}^d).
\end{equation}
Then its range space
\begin{equation}
V_{K,p}=\Big\{Tf:f\in L_\nu^p(\mathbb{R}^d)\Big\}=\Big\{f\in L_\nu^p(\mathbb{R}^d):Tf=f\Big\}
\end{equation}
is a weighted reproducing kernel subspace of $ L_\nu^p(\mathbb{R}^d)$ \cite{J1,N1,X1}, which means that for any $x\in \mathbb{R}^d$, there exists a $C_x>0$ such that
\begin{equation}
|f(x)|\leq C_x\|f\|_{L_\nu^p(\mathbb{R}^d)},\ f\in V_{K,p}.
\end{equation}

Let $0<\delta<1$.  Define a  subset of $V_{K,p}$ by
\begin{equation}
V_{K,p}(R,\delta)=\bigg\{f\in V_{K,p}:\int_{C_{R}}|f(x)\nu(x)|^{p}dx\geq(1-\delta)\int_{\mathbb{R}^{d}}|f(x)\nu(x)|^{p}dx\bigg\},
\end{equation}
which contains all functions in $V_{K,p}$ whose energy concentrate on the cube $C_{R}$.

\par This paper is organized as follows. In section 2, we show that a function $f\in V_{K,p}$ can be approximated by a function $f_N$ in a finite dimensional subspace $V^{N}_{K,p}$ on any bounded domains.  In section 3, we give an estimate for the  covering number of normalized $V^{N}_{K,p}$. In section 4,  we  prove that the sampling inequality holds with high probability for all functions in $V_{K,p}(R,\delta)$.
%In section 5,  a reconstruction algorithm based on random samples is provided for functions in $V^{N}_{K,p}$.

\section{Approximation to $V_{K,p}$ }

\ \ \ \ In this section, we will show that  $V_{K,p}$ can be approximated by a finite dimensional subspace on any bounded domains.
The following definitions of frame is similar to \cite{A1,H2,S4}.
\begin{Df}\label{def1}
Let $V$ be a Banach subspace of $L_\nu^{p}(\mathbb{R}^d)$ and $1/p+1/p'=1$. A family $\Psi=\{\psi_\gamma\}_{\gamma\in \Gamma}$ of functions in
$L_{1/\nu}^{p'}(\mathbb{R}^d)$
 is a $p$-frame for $V$, if there exist positive constants $A_p$ and $B_p$ such that
\begin{equation*}
A_p\|f\|_{L_\nu^p}\leq \big\|\{\langle f,\psi_\gamma\rangle\}_{\gamma\in \Gamma}\big\|_{\ell_\nu^p}\leq B_p\|f\|_{L_\nu^p},\ \forall \ f\in V.
\end{equation*}
\end{Df}
\begin{Df}\label{def2}
Let  $V\subset L_\nu^p(\mathbb{R}^{d})$ and $W\subset L_{1/\nu}^{p'}(\mathbb{R}^{d})$. The $p$-frame $\widetilde{\Phi}=\{\widetilde{\phi}_\gamma\}_{\gamma\in \Gamma}\subset W$ for $V$ and the $p'$-frame $\Phi=\{\phi_\gamma\}_{\gamma\in \Gamma}\subset V$ for $W$ form a dual pair if the following reconstruction formulae hold:
\begin{equation}\label{r23}
f=\sum \limits_{\gamma\in \Gamma}\langle f,\widetilde{\phi}_\gamma\rangle\phi_\gamma \ for \ all\ f\in V
\end{equation}
and
\begin{equation}\label{r24}
g=\sum \limits_{\gamma\in \Gamma}\langle g,\phi_\gamma\rangle\widetilde{\phi}_\gamma \ for \ all\ g\in W.
\end{equation}
\end{Df}

\begin{Lm}\cite{N1}\label{lm2.3}
Let $1\leq p\leq \infty$, $T$ be an idempotent integral operator on $L_\nu^p(\mathbb{R}^d)$ whose kernel $K$ satisfies \eqref{lab1} and \eqref{lab2}, and let $V_{K,p}$
be the range space of $T$. Then there exists a relatively-separated subset $\Lambda=\delta_0\mathbb{Z}^d$ with $\delta_0$ being determined by the condition $\|K\|_{\mathcal{W}}\|\omega_{\delta_0}(K)\|_{\mathcal{W}}<1$, and two families $\Phi=\{\phi_\lambda\}_{\lambda\in \Lambda}$ in $V_{K,p}$ and $\widetilde{\Phi}=\{\widetilde{\phi}_\lambda\}_{\lambda\in \Lambda}$ in $V_{K,p}^\ast$ which are defined by
\begin{equation}\label{3.26}
\phi_\lambda(x)=\delta_0^{-d/p}\int_{\mathbb{R}^d}\int_{[-\delta_0/2,\delta_0/2]^d}K_{\delta_0}(x,z_1)K(z_1,\lambda+z_2)dz_2dz_1
\end{equation}
with \begin{equation}
K_{\delta_0}(x,y)=\delta_0^{-d}\int_{[-\delta_0/2,\delta_0/2]^d}\int_{[-\delta_0/2,\delta_0/2]^d}\sum \limits_{\lambda\in \delta_0 \mathbb{Z}^d}
K(x,\lambda+z_1)K(\lambda+z_2,y)dz_1dz_2,
\end{equation}
and \begin{equation}
\widetilde{\phi}_\lambda(x)=\delta_0^{-d+d/p}\int_{[-\delta_0/2,\delta_0/2]^d}K(\lambda+z,x)dz
\end{equation}
such that
\begin{enumerate}
\item[(i)] Both $\Phi$ and $\widetilde{\Phi}$ are localized in the sense that
\begin{equation}\label{j4}
|\phi_\lambda(x)|+|\widetilde{\phi}_\lambda(x)|\leq h(x-\lambda),
\end{equation}
where $h\in L_\omega ^1(\mathbb{R}^d)$.
\item[(ii)]  $\Phi$ and $\widetilde{\Phi}$ form a dual frame pair for $V_{K,p}$ and $V_{K,p}^\ast$.
\item[(iii)]  Both $V_{K,p}$ and $V_{K,p}^\ast$ are generated by  $\Phi$ and $\widetilde{\Phi}$ in the sense that
\begin{equation}\label{3.25}
V_{K,p}=\Big\{\sum \limits_{\lambda\in \Lambda}c(\lambda)\phi_\lambda:\ (c(\lambda))_{\lambda\in \Lambda}\in \ell_\nu^p(\Lambda)\Big\}
\end{equation}
and
\begin{equation}
V_{K,p}^\ast=\Big\{\sum \limits_{\lambda\in \Lambda}\widetilde{c}(\lambda)\widetilde{\phi}_\lambda:\ (\widetilde{c}(\lambda))_{\lambda\in \Lambda}\in \ell_{1/\nu}^{p/(p-1)}(\Lambda)\Big\}.
\end{equation}
\item[(iv)] $\|K_{\delta_0}\|_{\mathcal{W}}<\infty$ and $\lim \limits_{\delta\rightarrow 0}\|\omega_\delta(K_{\delta_0})\|_{\mathcal{W}}=0$.
\end{enumerate}
\end{Lm}

Based on Lemma \ref{lm2.3}, for a given positive integer $N$, define a finite dimensional subspace
\begin{equation}\label{31}\
V_{K,p}^N=\bigg\{\sum\limits_{\lambda\in \Lambda\cap [-N,N]^d}c(\lambda)\phi_{\lambda}: c(\lambda)\in \mathbb{R}\bigg\}
\end{equation} of $V_{K,p}$
and its normalization
\begin{equation}\label{32}
V_{K,p}^{N,\ast}=\bigg\{f\in V_{K,p}^N:\|f\|_{L_\nu^p(\mathbb{R}^d)}=1\bigg\}.
\end{equation}

In the following, we will show that $V_{K,p}$ can be approximated by $ V_{K,p}^N$ on any bounded domains $C_M=[-M,M]^d$ with $M>0$.
%\begin{Lm}
%Let $\delta_0$ be a positive number satisfying $\delta_0<1/d$. If $K(x,y)$ satisfies the condition \eqref{jz1}, then one has
%\begin{equation}
%\lim \limits_{N\rightarrow \infty}\sup \limits_{x\in C_M} \sum \limits_{\lambda=\delta_0k,|k|>N}\int_{[-\delta_0/2,\delta_0/2]^d}|K(x,\lambda+y)|dy=0.
%\end{equation}
%\end{Lm}

\begin{Lm}\label{lm2.2}
Let $1\leq p\leq \infty$ and $p'$ be the conjugate number of $p$. Suppose that $K$ satisfies the assumptions \eqref{jz1} and \eqref{lab2}. If $f\in V_{K,p}$ and $\|f\|_{L_\nu^p(\mathbb{R}^d)}=1$, then for any given $\varepsilon>0$, there exist $N=N(\varepsilon,M)$ and
$f_{N}\in V_{K,p}^{N}$ such that
\begin{equation}\label{eq3}
\|f-f_{N}\|_{L_\nu^p(C_{M})}\leq\varepsilon \ and\ \|f-f_{N}\|_{L_\nu^{\infty}(C_{M})}\leq\frac{\varepsilon}{(2M)^{d/p}}.
\end{equation}
\end{Lm}
{\bf Proof} \ Since $f\in V_{K,p}$, it follows from \eqref{3.25} that  $f=\sum \limits_{\lambda\in \Lambda}\langle f, \widetilde{\phi}_\lambda\rangle\phi_\lambda$ for $\Lambda=\delta_0\mathbb{Z}^d$ with $\delta_0$ being chosen such that $\|K\|_{\mathcal{W}}\|\omega_{\delta_0}(K)\|_{\mathcal{W}}<1$. Take
\begin{equation}\label{hc1}
f_N=\sum\limits_{\lambda\in \Lambda\cap [-N,N]^d}\langle f, \widetilde{\phi}_\lambda\rangle\phi_{\lambda}\in V_{K,p}^{N}.
\end{equation}
%Since $\|K\|_{\mathcal{W}}<\infty$  and $\|K_{\delta_0}\|_{\mathcal{W}}<\infty$, then
%\begin{align}\label{j2}
%&\ \ \ \ \sum \limits_{\lambda\in \Lambda}|\phi_\lambda(x)|\omega(x-\lambda)\nonumber\\
%&\leq \delta_0^{-d/p}\int_{\mathbb{R}^d}|K_{\delta_0}(x,z_1)|\omega(x-z_1)\sum \limits_{k\in \mathbb{Z}^d}\int_{\delta_0k+[-\frac{\delta_0}{2},\frac{\delta_0}{2}]^d}
%|K(z_1,z_2)|\omega(z_1-\delta_0k)dz_2dz_1\nonumber\\
%&\leq \delta_0^{-d/p}\Big(\max \limits_{x\in[-\frac{\delta_0}{2},\frac{\delta_0}{2}]^d}\omega(x)\Big)\|K_{\delta_0}\|_{\mathcal{W}}\|K\|_{\mathcal{W}}.
%\end{align}
 For $k=(k_1,k_2,\cdots,k_d)\in \mathbb{Z}^d$, let $|k|=\max\{|k_1|,|k_2|,\cdots,|k_d|\}$. Then
\begin{align}\label{eq1}
&\ \ \ \ \ |f(x)-f_N(x)|\nu(x)\nonumber\\
%&=\Big|\sum \limits_{\lambda\in \Lambda\cap\{\mathbb{R}^d\backslash [-N,N]^d\}}\langle f,\widetilde{\phi}_\lambda\rangle\phi_\lambda(x)\nu(x)\Big|\nonumber\\
&\leq \sum \limits_{\lambda\in \Lambda\cap\{\mathbb{R}^d\backslash [-N,N]^d\}}|\langle f,\widetilde{\phi}_\lambda\rangle|\cdot|\phi_\lambda(x)|\nu(x)\nonumber\\
&\leq C_0\sum \limits_{\lambda\in \Lambda\cap\{\mathbb{R}^d\backslash [-N,N]^d\}}|\langle f,\widetilde{\phi}_\lambda\rangle|\nu(\lambda)\cdot|\phi_\lambda(x)|\omega(x-\lambda)\nonumber\\
&\leq C_0\Big\|\big\{\langle f,\widetilde{\phi}_\lambda \rangle \nu(\lambda)\big\}_{\lambda\in \Lambda\cap\{\mathbb{R}^d\backslash [-N,N]^d\}}\Big\|_{\ell^{p}}\Big\|\big\{\phi_\lambda(x)\omega(x-\lambda)\big\}_{\lambda\in \Lambda\cap\{\mathbb{R}^d\backslash [-N,N]^d\}}\Big\|_{\ell^{p'}}\nonumber\\
&\leq C_0\Big\|\big\{\langle f,\widetilde{\phi}_\lambda\rangle\big\}_{\lambda\in \Lambda}\Big\|_{\ell_\nu^p}\Big(\sum \limits_{\lambda\in \Lambda\cap\{\mathbb{R}^d\backslash [-N,N]^d\}}|\phi_\lambda(x)|\omega(x-\lambda)\Big).
\end{align}
Since $\widetilde{\Phi}=\{\widetilde{\phi}_\lambda\}_{\lambda\in \Lambda}$ is a $p$-frame of $V_{K,p}$,
by Definition \ref{def1}, one has
\begin{equation}
\big\|\big\{\langle f,\widetilde{\phi}_\lambda\rangle\big\}_{\lambda\in \Lambda}\big\|_{\ell_\nu^p}\leq B_p\|f\|_{L_\nu^p(\mathbb{R}^d)}=B_p.
\end{equation}
Moreover, it follows from \eqref{3.26} that
\begin{align}\label{3.26.1}
& \ \ \ \ \sum \limits_{\lambda\in \Lambda\cap\{\mathbb{R}^d\backslash [-N,N]^d\}}|\phi_\lambda(x)|\omega(x-\lambda)\nonumber\\
%&\leq \delta_0^{-d/p}\sum \limits_{\lambda\in \Lambda\cap\{\mathbb{R}^d\backslash [-N,N]^d\}}\int_{\mathbb{R}^d}\int_{[-\delta_0/2,\delta_0/2]^d}|K_{\delta_0}(x,z_1)||K(z_1,\lambda+z_2)|\omega(x-\lambda)dz_2dz_1\nonumber\\
&\leq \delta_0^{-d/p}\int_{\mathbb{R}^d}|K_{\delta_0}(x,z_1)|\omega(x-z_1)\sum \limits_{\lambda=\delta_0k,|k|>N}\int_{[-\delta_0/2,\delta_0/2]^d+\lambda}|K(z_1,z_2)|\omega(z_1-\lambda)dz_2dz_1\nonumber\\
&\leq \delta_0^{-d/p}\big(\mathop {\max }\limits_{x
 \in [-\delta_0/2,\delta_0/2]^d}\omega(x )\big)\int_{\mathbb{R}^d}|K_{\delta_0}(x,z_1)|\omega(x-z_1)\cdot\nonumber\\
 & \ \ \ \ \sum \limits_{\lambda=\delta_0k,|k|>N}\int_{[-\delta_0/2,\delta_0/2]^d+\lambda}|K(z_1,z_2)|\omega(z_2-z_1)dz_2dz_1\nonumber\\
&\leq \delta_0^{-d/p}\big(\mathop {\max }\limits_{x
 \in [-\delta_0/2,\delta_0/2]^d}\omega(x )\big)\int_{\mathbb{R}^d}|K_{\delta_0}(x,z_1)|\omega(x-z_1)\cdot\nonumber\\
 &\ \ \ \ \sum \limits_{\lambda=\delta_0k,|k|>N}\int_{[-\delta_0/2,\delta_0/2]^d+\lambda}\frac{\widetilde{C}}{(1+\|z_{1}\|_1+\|z_2\|_1)^\alpha}dz_2dz_1\nonumber\\
&\leq \delta_0^{-d/p}\widetilde{C}\big(\mathop {\max }\limits_{x
 \in [-\delta_0/2,\delta_0/2]^d}\omega(x )\big)\|K_{\delta_0}\|_{\mathcal{W}}\sum \limits_{\lambda=\delta_0k,|k|>N}\int_{[-\delta_0/2,\delta_0/2]^d+\lambda}\frac{1}{(1+\|z_2\|_1)^\alpha}dz_2.
 \end{align}
Since $\alpha>d$ and $ \lim \limits_{N\rightarrow\infty}\sum \limits_{\lambda=\delta_0k,|k|>N}\int_{[-\delta_0/2,\delta_0/2]^d+\lambda}\frac{1}{(1+\|z_2\|_1)^\alpha}dz_2=0$ is independent of the variable $x$, this together with \eqref{eq1}-\eqref{3.26.1} obtains the desired result.

\begin{Lm}\label{lm2.5}
Suppose that $K$ satisfies the assumptions \eqref{lab1} and \eqref{lab2}, then there exists a positive constant $C_K=C_0\Big[\delta_0^{-d/p}\big(\mathop {\max }\limits_{x \in [-\delta_0/2,\delta_0/2]^d}\omega(x)\big)\|K_{\delta_0}\|_{\mathcal{W}}\|K\|_{\mathcal{W}}\Big]^{1-1/p}\|h\|_{L_\omega^1}^{1/p}$ such that
\begin{equation}
\Big\|\sum \limits_{\lambda\in \Lambda}c(\lambda)\phi_\lambda\Big\|_{L_\nu^p(\mathbb{R}^d)}\leq C_K\big\|\big(c(\lambda)\big)_{\lambda\in \Lambda}\big\|_{\ell_\nu^p(\Lambda)}.
\end{equation}
\end{Lm}
{\bf Proof}  It follows from  \eqref{3.26} that
\begin{align}\label{3.26.2}
& \ \ \ \ \sum \limits_{\lambda\in \Lambda}|\phi_\lambda(x)|\omega(x-\lambda)\nonumber\\
&\leq \delta_0^{-d/p}\int_{\mathbb{R}^d}|K_{\delta_0}(x,z_1)|\omega(x-z_1)\sum \limits_{\lambda=\delta_0k}\int_{[-\delta_0/2,\delta_0/2]^d+\lambda}|K(z_1,z_2)|\omega(z_1-\lambda)dz_2dz_1\nonumber\\
&\leq \delta_0^{-d/p}\big(\mathop {\max }\limits_{x \in [-\delta_0/2,\delta_0/2]^d}\omega(x)\big)\|K_{\delta_0}\|_{\mathcal{W}}\|K\|_{\mathcal{W}}.
\end{align}
If $1\leq p<\infty$, by  \eqref{j4} and \eqref{3.26.2}, we can obtain
\begin{align*}
  & \ \ \ \ \ \ \Big\|\sum \limits_{\lambda\in \Lambda}c(\lambda)\phi_\lambda\Big\|_{L_\nu^p(\mathbb{R}^d)}^p\nonumber\\
  &\leq \int_{\mathbb{R}^d}\Big(\sum \limits_{\lambda\in \Lambda}|c(\lambda)|\cdot|\phi_\lambda(x)|\nu(x)\Big)^pdx\nonumber\\
&\leq C_0^p\int_{\mathbb{R}^d}\Big(\sum \limits_{\lambda\in \Lambda}|c(\lambda)|\nu(\lambda)\cdot|\phi_\lambda(x)|\omega(x-\lambda)\Big)^pdx\nonumber\\
&\leq C_0^p\int_{\mathbb{R}^d}\Big(\sum \limits_{\lambda\in \Lambda}|c(\lambda)|^p\nu(\lambda)^p\cdot|\phi_\lambda(x)|\omega(x-\lambda)\Big)\Big(\sum \limits_{\lambda\in \Lambda}|\phi_\lambda(x)|\omega(x-\lambda)\Big)^{p/p'}dx\nonumber\\
&\leq C_0^p\Big[\delta_0^{-d/p}\big(\mathop {\max }\limits_{x \in [-\delta_0/2,\delta_0/2]^d}\omega(x)\big)\|K_{\delta_0}\|_{\mathcal{W}}\|K\|_{\mathcal{W}}\Big]^{p-1}\sum \limits_{\lambda\in \Lambda}|c(\lambda)|^p\nu(\lambda)^p\int_{\mathbb{R}^d}|\phi_\lambda(x)|\omega(x-\lambda)dx\nonumber\\
&\leq C_0^p\Big[\delta_0^{-d/p}\big(\mathop {\max }\limits_{x \in [-\delta_0/2,\delta_0/2]^d}\omega(x)\big)\|K_{\delta_0}\|_{\mathcal{W}}\|K\|_{\mathcal{W}}\Big]^{p-1}\|h\|_{L_\omega ^1}\big\|\big(c(\lambda)\big)_{\lambda\in \Lambda}\big\|_{\ell_\nu^p(\Lambda)}^p\nonumber\\
&=C_K^p\big\|\big(c(\lambda)\big)_{\lambda\in \Lambda}\big\|_{\ell_\nu^p(\Lambda)}^p.
\end{align*}
If $p=\infty$, then it follows from \eqref{3.26.2} that
\begin{equation*}
 \Big\|\sum \limits_{\lambda\in \Lambda}c(\lambda)\phi_\lambda\Big\|_{L_\nu^{\infty}(\mathbb{R}^d)}\leq C_0\Big(\sum \limits_{\lambda\in \Lambda}|\phi_\lambda(x)|\omega(x-\lambda)\Big)\big\|\big(c(\lambda)\big)_{\lambda\in \Lambda}\big\|_{\ell_\nu^\infty(\Lambda)}\leq C_K\big\|\big(c(\lambda)\big)_{\lambda\in \Lambda}\big\|_{\ell_\nu^\infty(\Lambda)}.
\end{equation*}
\section{Covering number for $V_{K,p}^{N,\ast}$}
\ \ \ \ \ In this section, we discuss the covering number of $V_{K,p}^{N,\ast}$ with respect to the norm $\|\cdot\|_{L_\nu ^{\infty}(\mathbb{R}^d)}$.
Let $S$ be a metric space and $\eta>0$, the covering number $\mathcal{N}(S,\eta)$ is defined to be the minimal integer $m\in \mathbb{N}$ such that there exist $m$ disks with radius $\eta$ covering $S$.
%\subsection{Covering number of ${\bf \widetilde{V}_{N,K}(\Phi)}$}
\begin{Lm}\label{Lm 3.1}(\cite{C1})
Suppose $\mathbb{E}$ is a finite dimensional Banach space with dim$\mathbb{E}=s$. Let $B_{\varepsilon}:=\{x\in \mathbb{E}:\|x\|\leq \varepsilon\}$ be the closed ball of radius $\varepsilon$ centered at the origin. Then
\begin{equation*}
\mathcal{N}(B_{\varepsilon},\eta)\leq \Big(\frac{2\varepsilon}{\eta}+1\Big)^{s}.
\end{equation*}
\end{Lm}
\par Note that
\begin{equation}
dim\big(V_{K,p}^{N}\big)\leq \sharp\big\{\lambda\in \Lambda:\lambda\in[-N,N]^d\big\}\leq \big(\frac{2N}{\delta_0}+1\big)^d.
\end{equation}
Then by Lemma \ref{Lm 3.1}, we have the following result.
\begin{Lm}\label{Lm 3.5}
Let $V_{K,p}^{N,\ast}$ be defined by \eqref{32}. Then for any $\eta>0$, the covering number of $V_{K,p}^{N,\ast}$ concerning the norm $\|\cdot\|_{L_\nu^p(\mathbb{R}^d)}$ is bounded by
\begin{equation*}
\mathcal{N}\big(V_{K,p}^{N,\ast},\eta\big)\leq\exp\bigg(\big(\frac{2N}{\delta_0}+1\big)^d\ln\Big(\frac{2}{\eta}+1\Big)\bigg).
\end{equation*}
\end{Lm}
\begin{Lm}\label{Lm 3.3}
Suppose that  $K$ satisfies the assumptions \eqref{lab1} and \eqref{lab2}. Then for every $f\in V_{K,p}$, we have
\begin{equation}\label{eq9}
\|f\|_{L_\nu^{\infty}(\mathbb{R}^d)}\leq C^\ast\|f\|_{L_\nu^p(\mathbb{R}^d)},
\end{equation}
where
\begin{equation}
C^\ast= B_pC_0\delta_0^{-d/p}\big(\mathop {\max }\limits_{x \in [-\delta_0/2,\delta_0/2]^d}\omega(x)\big)
\|K_{\delta_0}\|_{\mathcal{W}}\|K\|_{\mathcal{W}}.
\end{equation}
\end{Lm}
{\bf Proof}\  Suppose that
$f\in V_{K,p}$, then it follows from Definition \ref{def1}, Definition \ref{def2} and Lemma \ref{lm2.3} that $f=\sum \limits_{\lambda\in \Lambda}\langle f, \widetilde{\phi}_\lambda\rangle\phi_\lambda$. Moreover, we can obtain from \eqref{3.26.2}  that
\begin{align*}\label{eq7}
\|f\|_{L_\nu^{\infty}(\mathbb{R}^d)}&\leq\sup\limits_{x\in \mathbb{R}^{d}}\sum\limits_{\lambda\in \Lambda}|\langle f, \widetilde{\phi}_\lambda\rangle||\phi_\lambda(x)|\nu(x)\nonumber\\
&\leq C_0\sup\limits_{x\in \mathbb{R}^{d}}\sum\limits_{\lambda\in \Lambda}|\langle f, \widetilde{\phi}_\lambda\rangle|\nu(\lambda)\cdot|\phi_\lambda(x)|\omega(x-\lambda)\nonumber\\
&\leq C_0\delta_0^{-d/p}\big(\mathop {\max }\limits_{x \in [-\delta_0/2,\delta_0/2]^d}\omega(x)\big)\|K_{\delta_0}\|_{\mathcal{W}}\|K\|_{\mathcal{W}}\big\|\big\{\langle f,\widetilde{\phi}_\lambda\rangle\big\}_{\lambda\in \Lambda}\big\|_{\ell_\nu^p}\nonumber\\
&\leq B_pC_0\delta_0^{-d/p}\big(\mathop {\max }\limits_{x \in [-\delta_0/2,\delta_0/2]^d}\omega(x)\big)\|K_{\delta_0}\|_{\mathcal{W}}\|K\|_{\mathcal{W}}\|f\|_{L_\nu^{p}(\mathbb{R}^d)}.
\end{align*}

\begin{Lm}\label{lem 3.7}
Suppose that $K$ satisfies the assumptions \eqref{lab1} and \eqref{lab2},  then the covering number of $V_{K,p}^{N,\ast}$ with respect to $\|\cdot\|_{L_\nu^{\infty}(\mathbb{R}^d)}$ is bounded by
\begin{equation*}
\mathcal{N}\big(V_{K,p}^{N,\ast},\eta\big)\leq\exp\bigg(\big(\frac{2N}{\delta_0}+1\big)^d\ln\Big(\frac{2C^{*}}{\eta}+1\Big)\bigg).
\end{equation*}
\end{Lm}
{\bf Proof} \  By Lemma \ref{Lm 3.5}, the covering number of $V_{K,p}^{N,\ast}$ with respect to $\|\cdot\|_{L_\nu^{p}(\mathbb{R}^d)}$ satisfies
\begin{equation}\label{320}
\mathcal{N}\Big(V_{K,p}^{N,\ast},\frac{\eta}{C^{*}}\Big)\leq\exp\bigg(\big(\frac{2N}{\delta_0}+1\big)^d\ln\Big(\frac{2C^{*}}{\eta}+1\Big)\bigg).
\end{equation}
Let $\mathcal{F}$ be the corresponding $\frac{\eta}{C^{*}}$-net for $V_{K,p}^{N,\ast}$. It means that for every $f\in V_{K,p}^{N,\ast}$, there exists a $\widetilde{f}\in \mathcal{F}$ such that $\|f-\widetilde{f}\|_{L_\nu^{p}(\mathbb{R}^d)}\leq \frac{\eta}{C^{*}}$. By Lemma \ref{Lm 3.3}, we have
\begin{equation*}
\|f-\widetilde{f}\|_{L_\nu^{\infty}(\mathbb{R}^d)}\leq C^{*}\|f-\widetilde{f}\|_{L_\nu^{p}(\mathbb{R}^d)}\leq\eta.
\end{equation*}
Therefore, $\mathcal{F}$ is also a $\eta$-net of $V_{K,p}^{N,\ast}$ with respect to the norm $\|\cdot\|_{L_\nu^{\infty}(\mathbb{R}^d)}$. Since
\begin{equation*}
\sharp(\mathcal{F})\leq \exp\bigg(\big(\frac{2N}{\delta_0}+1\big)^d\ln\Big(\frac{2C^{*}}{\eta}+1\Big)\bigg),
\end{equation*}
the desired result is proved.

\section{Random sampling inequality of $V_{K,p}(R,\delta)$}
\ \ \ \ \ Let $X=\{x_{j}:j\in \mathbb{N}\}$ be a sequence of independent random  variables that are drawn from a general probability distribution over $\mathbb{R}^d$ with density function $\rho$ satisfying \eqref{den}.
 Then for any $f\in V_{K,p}$, we introduce the random variables
\begin{equation}\label{b1}
X_{j}(f)=|f(x_{j})\nu(x_{j})|^{p}-\int_{\mathbb{R}^d}\rho(x)|f(x)\nu(x)|^{p}dx.
\end{equation}
 It is easy to see that $X_{j}(f)$ is a sequence of independent random variables with expectation $\mathbb{E}[X_{j}(f)]=0$. Next, we will give some estimates for $X_{j}(f)$.
\begin{Lm}\label{Lm 3.8}
Let $\rho(x)$ be a probability density function over $\mathbb{R}^d$ satisfying \eqref{den}. Then for any $f,g\in V_{K,p}$, the following inequalities hold:\\
$(1)\ \|X_{j}(f)\|_{\ell^{\infty}}\leq \|f\|_{L_\nu^{\infty}(\mathbb{R}^d)}^{p}$.\\
$(2)\ \|X_{j}(f)-X_{j}(g)\|_{\ell^{\infty}}\leq 2p\Big(\max\big\{\|f\|_{L_\nu^{\infty}(\mathbb{R}^d)},\|g\|_{L_\nu^{\infty}(\mathbb{R}^d)}\big\}\Big)^{p-1}\|f-g\|_{L_\nu^{\infty}(\mathbb{R}^d)}$.\\
$(3)\ Var(X_{j}(f))\leq C_{\rho}\|f\|_{L_\nu^{\infty}(\mathbb{R}^d)}^{p}\|f\|_{L_\nu^{p}(\mathbb{R}^d)}^{p}$.\\
$(4)\ Var\big(X_{j}(f)-X_{j}(g)\big)\leq pC_{\rho}\Big(\max\big\{\|f\|_{L_\nu^{\infty}(\mathbb{R}^d)},\|g\|_{L_\nu^{\infty}(\mathbb{R}^d)}\big\}\Big)^{p-1}\|f-g\|_{L_\nu^{\infty}(\mathbb{R}^d)}
\Big(\|f\|_{L_\nu^{p}(\mathbb{R}^d)}^{p}+\|g\|_{L_\nu^{p}(\mathbb{R}^d)}^{p}\Big)$.
\end{Lm}
{\bf Proof}
(1) Direct computation obtains \begin{equation*}
\|X_{j}(f)\|_{\ell^{\infty}}
\leq\sup\limits_{x\in \mathbb{R}^d}\max\bigg\{
|f(x)\nu(x)|^{p},\int_{\mathbb{R}^d}\rho(x)|f(x)\nu(x)|^{p}dx\bigg\}
\leq\|f\|_{L_\nu^{\infty}(\mathbb{R}^d)}^{p}.
\end{equation*}
(2) By mean value theorem, one has
\begin{eqnarray*}
\|X_{j}(f)-X_{j}(g)\|_{\ell^{\infty}}
&\leq&\sup\limits_{x\in \mathbb{R}^d}\bigg(\Big||f(x)\nu(x)|^{p}-|g(x)\nu(x)|^{p}\Big|+\int_{\mathbb{R}^d}\rho(x)\Big||f(x)\nu(x)|^{p}-|g(x)\nu(x)|^{p}\Big|dx\bigg)\\
&\leq&2\sup\limits_{x\in \mathbb{R}^d}\Big||f(x)\nu(x)|^{p}-|g(x)\nu(x)|^{p}\Big|\\
&=&2p\Big(\max\big\{\|f\|_{L_\nu^{\infty}(\mathbb{R}^d)},\|g\|_{L_\nu^{\infty}(\mathbb{R}^d)}\big\}\Big)^{p-1}\|f-g\|_{L_\nu^{\infty}(\mathbb{R}^d)}.
\end{eqnarray*}
(3) Since $\mathbb{E}[X_{j}(f)]=0$, then
\begin{eqnarray*}
Var(X_{j}(f))&=&\mathbb{E}[(X_{j}(f))^{2}]\\
&=&\mathbb{E}\big[|f(x_{j})\nu(x_{j})|^{2p}\big]-\Big(\int_{\mathbb{R}^d}\rho(x)|f(x)\nu(x)|^{p}dx\Big)^{2}\\
&\leq&\int_{\mathbb{R}^d}\rho(x)|f(x)\nu(x)|^{2p}dx\\
&\leq&C_{\rho}\|f\|_{L_\nu^{\infty}(\mathbb{R}^d)}^{p}\|f\|_{L_\nu^{p}(\mathbb{R}^d)}^{p}.
\end{eqnarray*}
(4) Using the similar method as $(3)$, we have
\begin{eqnarray*}
&&Var\big(X_{j}(f)-X_{j}(g)\big)\\&=&\mathbb{E}\big[\big(X_{j}(f)-X_{j}(g)\big)^{2}\big]\\
&\leq& C_{\rho}\int_{\mathbb{R}^d}\Big(|f(x)\nu(x)|^{p}-|g(x)\nu(x)|^{p}\Big)^{2}dx\\
&\leq& C_{\rho}\int_{\mathbb{R}^d}\Big||f(x)\nu(x)|^{p}-|g(x)\nu(x)|^{p}\Big|\Big(|f(x)\nu(x)|^{p}+|g(x)\nu(x)|^{p}\Big)dx\\
&\leq& C_{\rho}\sup\limits_{x\in \mathbb{R}^d}\Big||f(x)\nu(x)|^{p}-|g(x)\nu(x)|^{p}\Big|\Big(\|f\|_{L_\nu^{p}(\mathbb{R}^d)}^{p}+\|g\|_{L_\nu^{p}(\mathbb{R}^d)}^{p}\Big)\\
&\leq& pC_{\rho}\Big(\max\big\{\|f\|_{L_\nu^{\infty}(\mathbb{R}^d)},\|g\|_{L_\nu^{\infty}(\mathbb{R}^d)}\big\}\Big)^{p-1}\|f-g\|_{L_\nu^{\infty}(\mathbb{R}^d)}
\Big(\|f\|_{L_\nu^{p}(\mathbb{R}^d)}^{p}+\|g\|_{L_\nu^{p}(\mathbb{R}^d)}^{p}\Big).
\end{eqnarray*}
\par  In the following lemma, we will show that a uniform large deviation inequality holds for functions in $V^{N,*}_{K,p}$ by Bernstein's inequality.
\begin{Lm}(Bernstein's inequality)\label{Lm 3.9}(\cite{B3})
Let $X_{1},X_{2},\ldots, X_{n}$ be independent random variables with expected values $\mathbb{E}(X_{j})=0$ for $j=1,2,\ldots, n$. Assume that $Var(X_{j})\leq \sigma^{2}$ and $|X_{j}|\leq M_0$ almost surely for all $j$. Then for any $\lambda\geq 0$,
\begin{equation*}
Prob\bigg(\bigg|\sum\limits_{j=1}^{n}X_{j}\bigg|\geq\lambda\bigg)\leq2\exp\bigg(-\frac{\lambda^{2}}{2n\sigma^{2}+\frac{2}{3}M_0\lambda}\bigg).
\end{equation*}
\end{Lm}
\begin{Lm}\label{lem 3.10}
 Let $\{x_{j}:j\in \mathbb{N}\}$ be a sequence of independent random variables  that are drawn from a general probability distribution over $\mathbb{R}^d$ with
 density function $\rho$ satisfying \eqref{den}. If $f\in V^{N,*}_{K,p}$, then for  $n\in \mathbb{N}$ and $\lambda\geq 0$,
\begin{equation*}
Prob\bigg(\sup\limits_{f\in V^{N,*}_{K,p}}\bigg|\sum\limits_{j=1}^{n}X_{j}(f)\bigg|\geq\lambda\bigg)\leq A\exp\bigg(-B\frac{\lambda^{2}}{12nC_{\rho}+2\lambda }\bigg),
\end{equation*}
where $A$ is of order $\exp(CN^{d})$ with  $C$ depending  on $\Lambda$ and $K$, and $B=\min\{\frac{\sqrt{2}}{2592p(C^{*})^{p-1}},\frac{3}{2(C^{*})^{p}}\}$.
\end{Lm}
{\bf Proof} \ For given $\ell\in \mathbb{N}$, we construct a $2^{-\ell}$-covering for $V^{N,*}_{K,p}$ with respect to the norm $\|\cdot\|_{L_\nu^{\infty}(\mathbb{R}^d)}$. Let $\mathcal{C}_{\ell}$ be the corresponding $2^{-\ell}$-net for $\ell=1,2,\ldots$. Then, $$\sharp(\mathcal{C}_{\ell})\leq\mathcal{N}\big(V^{N,*}_{K,p},2^{-\ell}\big).$$
\par For given $f\in V^{N,*}_{K,p}$, let $f_{\ell}$ be the function in $\mathcal{C}_{\ell}$ that is closest to $f$ with respect to the norm $\|\cdot\|_{L_\nu^{\infty}(\mathbb{R}^d)}$. Then, $\|f-f_{\ell}\|_{L_\nu^{\infty}(\mathbb{R}^d)}\leq 2^{-\ell}\rightarrow 0$ when $\ell\rightarrow \infty$. Moreover, by Lemma \ref{Lm 3.3} and the item $(2)$ of Lemma \ref{Lm 3.8}, we have
\begin{equation*}
X_{j}(f)=X_{j}(f_{1})+(X_{j}(f_{2})-X_{j}(f_{1}))+(X_{j}(f_{3})-X_{j}(f_{2}))+\cdots .
\end{equation*}
If $\sup\limits_{f\in V^{N,*}_{K,p}}\Big|\sum\limits_{j=1}^{n}X_{j}(f)\Big|\geq \lambda$, the event $\omega_{\ell}$ must hold for some $\ell\geq 1$, where
\begin{equation*}
\omega_{1}=\bigg\{\text{there exists} \ f_{1}\in\mathcal{C}_{1} \ \text{such that}\bigg |\sum\limits_{j=1}^{n}X_{j}(f_{1})\bigg|\geq \frac{\lambda}{2} \bigg\}
\end{equation*}
and for $\ell\geq 2$,
\begin{equation*}
\begin{aligned}
\omega_{\ell}=&\bigg\{\text{there exist}\  f_{\ell}\in \mathcal{C}_{\ell}\ \text{and} \ f_{\ell-1}\in\mathcal{C}_{\ell-1}\ \text{with}\ \|f_{\ell}-f_{\ell-1}\|_{L_\nu^{\infty}(\mathbb{R}^d)}\leq 3\cdot 2^{-\ell},\\&
\ \text{such that}\ \left|\sum_{j=1}^{n}\big(X_{j}(f_{\ell})-X_{j}(f_{\ell-1})\big)\right|\geq\frac{\lambda}{2\ell^{2}}\bigg\}.
\end{aligned}
\end{equation*}
If this is not the case, then with $f_{0}=0$, we have
\begin{equation*}
\bigg|\sum\limits_{j=1}^{n}X_{j}(f)\bigg|\leq\sum\limits_{\ell=1}^{\infty}\bigg|\sum\limits_{j=1}^{n}(X_{j}(f_{\ell})-X_{j}(f_{\ell-1}))\bigg|\leq\sum\limits_{\ell=1}^{\infty}\frac{\lambda}{2\ell^{2}}=\frac{\pi^{2}\lambda}{12}\leq\lambda.
\end{equation*}
Next, we estimate the probability of each $\omega_{\ell}$. By Lemma \ref{Lm 3.3}, Lemma \ref{Lm 3.8} and Lemma \ref{Lm 3.9}, for every fixed $f\in \mathcal{C}_{1}$,
\begin{eqnarray*}
Prob\bigg(\bigg|\sum\limits_{j=1}^{n}X_{j}(f)\bigg|\geq\frac{\lambda}{2}\bigg)&\leq&2\exp\bigg(-\frac{(\frac{\lambda}{2})^{2}}{2n\|Var(X_{j}(f))\|_{\ell^\infty}+\frac{2}{3}\|X_{j}(f)\|_{\ell^{\infty}}\cdot\frac{\lambda}{2}}\bigg)\\
&\leq&2\exp\bigg(-\frac{\lambda^{2}}{8nC_{\rho}(C^{*})^{p}+\frac{4}{3}\lambda (C^{*})^{p}}\bigg).
\end{eqnarray*}
By Lemma \ref{lem 3.7}, there are at most
\begin{equation*}
\mathcal{N}\Big(V^{N,*}_{K,p},\frac{1}{2}\Big)\leq\exp\bigg(\big(\frac{2N}{\delta_0}+1\big)^d\ln(4C^{*}+1)\bigg)
\end{equation*}
functions in $\mathcal{C}_{1}$. Thus, the probability of $\omega_{1}$ is bounded by
\begin{eqnarray}\label{w1}
Prob(\omega_1)&\leq&2\exp\bigg(\big(\frac{2N}{\delta_0}+1\big)^d\ln(4C^{*}+1)\bigg)\exp\bigg(-\frac{\lambda^{2}}{8nC_{\rho}(C^{*})^{p}+\frac{4}{3}\lambda (C^{*})^{p}}\bigg)\nonumber\\
&=&2\exp\bigg(\big(\frac{2N}{\delta_0}+1\big)^d\ln(4C^{*}+1)\bigg)\exp\bigg(-\frac{\lambda^{2}}{\frac{2}{3}(C^{*})^{p}(12nC_{\rho}+2\lambda )}\bigg).
\end{eqnarray}
For $\ell\geq 2$, we estimate the probability of $\omega_{\ell}$ in a similar way. For $f\in \mathcal{C}_{\ell}$, $g\in \mathcal{C}_{\ell-1}$ and $\|f-g\|_{L_\nu^{\infty}(\mathbb{R}^d)}\leq 3\cdot 2^{-\ell}$, it follows from Lemma \ref{Lm 3.3}, Lemma \ref{Lm 3.8} and Lemma \ref{Lm 3.9} that
\begin{eqnarray*}
&&Prob\bigg(\bigg|\sum\limits_{j=1}^{n}(X_{j}(f)-X_{j}(g))\bigg|\geq\frac{\lambda}{2\ell^{2}}\bigg)
\\&\leq&2\exp\bigg(-\frac{(\frac{\lambda}{2\ell^{2}})^{2}}{2n\|Var(X_{j}(f)-X_{j}(g))\|_{\ell^\infty}+\frac{2}{3}\|X_{j}(f)-X_{j}(g)\|_{\ell^{\infty}}\cdot\frac{\lambda}{2\ell^{2}}}\bigg)\\
&\leq&2\exp\bigg(-\frac{\upsilon 2^{\ell}}{\ell^{4}}\bigg),
\end{eqnarray*}
where $\upsilon=\frac{\lambda^{2}}{4p(C^{*})^{p-1}(12nC_{\rho}+2\lambda )}$.
 There are at most $\mathcal{N}\big(V^{N,*}_{K,p},2^{-\ell}\big)$ functions in $\mathcal{C}_{\ell}$ and $\mathcal{N}\big(V^{N,*}_{K,p},2^{-\ell+1}\big)$ functions in $\mathcal{C}_{\ell-1}$. Therefore, we have
\begin{eqnarray*}
Prob\Big(\bigcup \limits_{\ell=2}^{\infty}\omega_{\ell}\Big)
&\leq & \sum\limits_{\ell=2}^{\infty}\mathcal{N}\big(V^{N,*}_{K,p},2^{-\ell}\big)\mathcal{N}\big(V^{N,*}_{K,p},2^{-\ell+1}\big)2\exp\bigg(-\frac{\upsilon 2^{\ell}}{\ell^{4}}\bigg)\\
&\leq&2(2C^{\ast}+1)^{2\big(\frac{2N}{\delta_0}+1\big)^d}\sum\limits_{\ell=2}^{\infty}\exp\bigg((2\ln2)\big(\frac{2N}{\delta_0}+1\big)^d\ell-\frac{\upsilon 2^{\ell}}{\ell^{4}}\bigg)\\
&=:&C_{1}\sum\limits_{\ell=2}^{\infty}\exp\bigg(C_{2}\ell-\frac{\upsilon 2^{\ell}}{\ell^{4}}\bigg)\\
&=&C_{1}\sum\limits_{\ell=2}^{\infty}\exp\bigg(-\upsilon 2^{\frac{\ell}{2}}\bigg(\frac{2^{\frac{\ell}{2}}}{\ell^{4}}-\frac{C_{2}\ell }{2^{\frac{\ell}{2}}\upsilon}\bigg)\bigg),
\end{eqnarray*}
where $C_{1}=2(2C^{\ast}+1)^{2\big(\frac{2N}{\delta_0}+1\big)^d}$ and $C_{2}=(2\ln2)\big(\frac{2N}{\delta_0}+1\big)^d$.

Let $C_{3}=\min\limits _{\ell \geq 2} \frac{2^{\frac{\ell}{2}}}{\ell^{4}}=\frac{1}{324}$ and
$C_{4}=\max\limits _{\ell \geq 2} \frac{8p(C^{*})^{p-1}\ell\ln 2}{2^{\frac{\ell}{2}}}=6 \sqrt{2} p(C^{*})^{p-1}\ln 2 $. Then
\begin{eqnarray*}
\frac{2^{\frac{\ell}{2}}}{\ell^{4}}-\frac{C_{2}\ell }{2^{\frac{\ell}{2}} \upsilon}&=&\frac{2^{ \frac{\ell}{2}}}{\ell^{4}}-\frac{8\ell p(C^{*})^{p-1} \big(\frac{2N}{\delta_0}+1\big)^d(12nC_{\rho}+2\lambda )\ln2}{2^{\frac{\ell}{2}}\lambda^{2}}\\
&\geq&\frac{1}{324}-\frac{C_{4}\big(\frac{2N}{\delta_0}+1\big)^d(12nC_{\rho}+2\lambda )}{\lambda^{2}}.
\end{eqnarray*}
We first consider the case that
\begin{equation}\label{10.7.1}
\frac{1}{324}-\frac{C_{4}\big(\frac{2N}{\delta_0}+1\big)^d(12nC_{\rho}+2\lambda )}{\lambda^{2}}> \frac{1}{648}.
\end{equation}
Since  $\sum\limits_{\ell=2}^{\infty}e^{-pa^{\ell}}\leq \frac{e^{-ap}}{pa\ln a}$ for $p,a>0$(\cite{S2}), then
\begin{eqnarray*}
  Prob\bigg(\bigcup_{\ell=2}^{\infty}\omega_{\ell}\bigg) &\leq& \frac{C_{1}\exp\bigg(-\sqrt{2}\upsilon\bigg(\frac{1}{324}-\frac{C_{4}\big(\frac{2N}{\delta_0}+1\big)^d(12nC_{\rho}+2\lambda )}{\lambda^{2}}\bigg)\bigg)} {\sqrt{2}\ln\sqrt{2}\cdot\upsilon\bigg(\frac{1}{324}-\frac{C_{4}\big(\frac{2N}{\delta_0}+1\big)^d(12nC_{\rho}+2\lambda )}{\lambda^{2}}\bigg)} \\
 &=& \frac{2(2C^{\ast}+1)^{2\big(\frac{2N}{\delta_0}+1\big)^d}}{\sqrt{2}\ln\sqrt{2}\cdot\upsilon\bigg(\frac{1}{324}-\frac{C_{4}\big(\frac{2N}{\delta_0}+1\big)^d(12nC_{\rho}+2\lambda )}{\lambda^{2}}\bigg)} \\
 &&\times\exp\bigg(-\sqrt{2}\upsilon\bigg(\frac{1}{324}-\frac{C_{4}\big(\frac{2N}{\delta_0}+1\big)^d(12nC_{\rho}+2\lambda )}{\lambda^{2}}\bigg)\bigg).
\end{eqnarray*}
Under the condition \eqref{10.7.1}, we have
\begin{eqnarray*}
& &\sqrt{2}\ln\sqrt{2}\cdot\upsilon\bigg(\frac{1}{324}-\frac{C_{4}\big(\frac{2N}{\delta_0}+1\big)^d(12nC_{\rho}+2\lambda )}{\lambda^{2}}\bigg)\\
&\geq&\frac{\sqrt{2}\ln\sqrt{2}C_{4}\big(\frac{2N}{\delta_0}+1\big)^d}{4p(C^{*})^{p-1}}\\
&\geq&3\ln\sqrt{2}\ln2.
\end{eqnarray*}
This together with the probability of $\omega_{1}$ in \eqref{w1} obtains
\begin{equation*}
Prob\bigg(\sup\limits_{f\in V^{N,*}_{K,p}}\left|\sum\limits_{j=1}^{n}X_{j}(f)\right|\geq\lambda\bigg)\leq Prob\bigg(\bigcup\limits_{\ell=1}^{\infty}\omega_{\ell}\bigg)\leq A\exp\bigg(-B\frac{\lambda^{2}}{12nC_{\rho}+2\lambda }\bigg).
\end{equation*}
Here, $A$ is of order $\exp\big(CN^{d}\big)$ with $C=2^{d+1}\big(1+\frac{1}{\delta_0}\big)^d\ln(2C^{*}+1)$ and $B=\min\big\{\frac{\sqrt{2}}{2592p(C^{*})^{p-1}},\frac{3}{2(C^{*})^{p}}\big\}$.
Finally, we consider the case that
\begin{equation*}
\frac{1}{324}-\frac{C_{4}\big(\frac{2N}{\delta_0}+1\big)^d(12nC_{\rho}+2\lambda )}{\lambda^{2}}\leq\frac{1}{648}.
\end{equation*}
In this case, we can choose $C\geq 648 C_{4} B2^{d}\big(1+\frac{1}{\delta_0}\big)^d$ such that $A\exp\bigg(-B\frac{\lambda^{2}}{12nC_{\rho}+2\lambda }\bigg)\geq 1$. This completes the proof.
\begin{Lm}\label{lem 3.11}
 Let $X=\{x_{j}:j\in \mathbb{N}\}$ be a sequence of independent random variables  that are drawn from a general probability distribution over $\mathbb{R}^d$ with
 density function $\rho$ satisfying \eqref{den}. Then for any  $\gamma>0$, the  inequality
\begin{equation}\label{gs}
 nc_{\rho}\bigg(\|f\|_{L_\nu^{p}(C_{R})}^{p}-\gamma\|f\|_{L_\nu^{p}(\mathbb{R}^d)}^{p}\bigg)\leq\sum\limits_{j=1}^{n}|f(x_{j})\nu(x_{j})|^{p}
\leq  n\big(c_{\rho}\gamma+C_\rho\big)\|f\|_{L_\nu^{p}(\mathbb{R}^d)}^{p}
\end{equation}
holds for  function $f\in V_{K,p}^{N}$ with probability at least
\begin{equation*}
1-A\exp\bigg(-B\frac{\gamma^{2}n c_{\rho}^{2}}{12C_\rho+2\gamma c_{\rho}}\bigg),
\end{equation*}
where $A$ and $B$ are as in Lemma \ref{lem 3.10}.
\end{Lm}
{\bf Proof}  \ It is obvious that every $f\in V_{K,p}^{N}$ satisfies the inequality \eqref{gs} if and only if $f/\|f\|_{L_\nu^{p}(\mathbb{R}^d)}$ does. So we assume that $\|f\|_{L_\nu^{p}(\mathbb{R}^d)}=1$, then $f\in V_{K,p}^{N,\ast}$.  The event
\begin{equation*}
D=\bigg\{\sup\limits_{f\in V_{K,p}^{N,\ast}}\bigg|\sum\limits_{j=1}^{n}X_{j}(f)\bigg|>\gamma n c_{\rho}\bigg\}
\end{equation*}
is the complement of
\begin{equation*}
\begin{aligned}
\widetilde{D}=&\bigg\{n\int_{\mathbb{R}^d}\rho(x)|f(x)\nu(x)|^{p}dx-\gamma n c_{\rho}\leq\sum\limits_{j=1}^{n}|f(x_{j})\nu(x_{j})|^{p}\\
&\leq\gamma n c_{\rho}+n\int_{\mathbb{R}^d}\rho(x)|f(x)\nu(x)|^{p}dx, \quad \forall f\in V_{K,p}^{N,\ast}\bigg\}\\
\subseteq&\bigg\{ nc_{\rho}\bigg(\|f\|_{L_\nu^{p}(C_{R})}^{p}-\gamma\|f\|_{L_\nu^{p}(\mathbb{R}^d)}^{p}\bigg)\leq\sum\limits_{j=1}^{n}|f(x_{j})\nu(x_{j})|^{p}\\&
\leq n\big(c_{\rho}\gamma+C_\rho\big)\|f\|_{L_\nu^{p}(\mathbb{R}^d)}^{p},\quad \forall f\in V_{K,p}^{N}\bigg\}=\overline{D}.
\end{aligned}
\end{equation*}
Using Lemma \ref{lem 3.10}, the sampling inequality \eqref{gs} holds for all $f\in V_{K,p}^{N}$ with probability
\begin{equation*}
Prob(\overline{D})\geq Prob(\widetilde{D})=1-Prob(D)\geq 1-A\exp\bigg(-B\frac{\gamma^{2}n c_{\rho}^{2}}{12C_\rho+2\gamma c_{\rho}}\bigg).
\end{equation*}

\par In the following, we will show that if the sampling size is sufficiently large, the sampling inequality holds with overwhelming probability for  functions in $V_{K,p}(R,\delta)$.
\begin{Tm}\label{Tm 4.5}
Let $X=\{x_{j}:j\in \mathbb{N}\}$ be a sequence of independent random variables  that are drawn from a general probability distribution over $\mathbb{R}^d$ with
 density function $\rho$ satisfying \eqref{den}. Suppose that  $M>R$ is a constant such that $\{x_j: j=1,2,\cdots,n\}\subseteq C_M$, then for any $0<\varepsilon,\gamma<1$ which satisfy
\begin{equation}\label{25}
L(\varepsilon,\gamma)=:c_\rho\Big(1-\delta-p(1+\varepsilon)^{p-1}\varepsilon-\gamma\big(B_pC_K\big)^p\Big)-p\Big(C^\ast+\frac{\varepsilon}{(2M)^{d/p}}\Big)^{p-1}\frac{\varepsilon}{(2M)^{d/p}}>0,
\end{equation}
 the  sampling inequality
\begin{equation}\label{24}
n L(\varepsilon,\gamma)\|f\|_{L_\nu^{p}(\mathbb{R}^d)}^{p}\leq\sum\limits_{j=1}^{n}|f(x_{j})\nu(x_{j})|^{p}\leq nU(\varepsilon,\gamma)\|f\|_{L_\nu^{p}(\mathbb{R}^d)}^{p}
\end{equation}
holds  for  function $f\in V_{K,p}(R,\delta)$ with probability at least
\begin{equation*}
1-A\exp\bigg(-B\frac{\gamma^{2}n c_{\rho}^{2}}{12C_\rho+2\gamma c_{\rho}}\bigg).
\end{equation*}
Here,  $U(\varepsilon,\gamma)=(c_\rho\gamma+C_\rho)\big(B_pC_K\big)^p+p\Big(C^\ast+\frac{\varepsilon}{(2M)^{d/p}}\Big)^{p-1}\frac{\varepsilon}{(2M)^{d/p}}$, $A$ and $B$ are the constants in Lemma \ref{lem 3.10} corresponding to $N=N(\varepsilon,M)$ in Lemma \ref{lm2.2}.
\end{Tm}
{\bf Proof } It is obvious that every $f\in V_{K,p}(R,\delta)$ satisfies the inequality \eqref{24} if and only if $f/\|f\|_{L_\nu^{p}(\mathbb{R}^d)}$ does. Hence, we assume that $\|f\|_{L_\nu^{p}(\mathbb{R}^d)}=1$.

%For random variables $\{x_j: j=1,2,\cdots,n\}$, there exists a $M>R$ such that $\{x_j: j=1,2,\cdots,n\}\subset C_M$.
For $\varepsilon>0$ satisfying \eqref{25}, it follows from Lemma \ref{lm2.2} that there exist positive integer $N=N(\varepsilon,M)$
and $f_{N}\in V^{N}_{K,p}$ such that
\begin{equation}
\|f-f_{N}\|_{L_\nu^{p}(C_{R})}\leq \|f-f_{N}\|_{L_\nu^{p}(C_{M})}\leq \varepsilon \ \ {\rm{and}} \ \ \|f-f_{N}\|_{L_\nu^{\infty}(C_{M})}\leq \frac{\varepsilon}{(2M)^{d/p}}.
\end{equation}
This together with mean value theorem and Lemma \ref{Lm 3.3} obtains
\begin{equation}\label{eq11}
\Big|\|f\|^p_{L_\nu^p(C_R)}-\|f_N\|^p_{L_\nu^p(C_R)}\Big|\leq p(1+\varepsilon)^{p-1}\varepsilon
\end{equation}
and
\begin{align}\label{eq10}
&\ \ \ \ \Big||f(x_j)\nu(x_j)|^p-|f_N(x_j)\nu(x_j)|^p\Big|\nonumber\\
&\leq p\Big(\max\big\{|f(x_j)\nu(x_j)|,|f_N(x_j)\nu(x_j)|\big\}\Big)^{p-1}|f(x_j)-f_N(x_j)|\nu(x_j)\nonumber\\
&\leq p\Big(C^\ast+\frac{\varepsilon}{(2M)^{d/p}}\Big)^{p-1}\frac{\varepsilon}{(2M)^{d/p}}.
\end{align}
It follows from \eqref{eq10} that
\begin{equation*}
\sum\limits_{j=1}^{n}|f_N(x_j)\nu(x_j)|^p-np\Big(C^\ast+\frac{\varepsilon}{(2M)^{d/p}}\Big)^{p-1}\frac{\varepsilon}{(2M)^{d/p}}\leq \sum\limits_{j=1}^{n}|f(x_j)\nu(x_j)|^p
\end{equation*}
\begin{equation}\label{eq12}
\leq \sum\limits_{j=1}^{n}|f_N(x_j)\nu(x_j)|^p+np\Big(C^\ast+\frac{\varepsilon}{(2M)^{d/p}}\Big)^{p-1}\frac{\varepsilon}{(2M)^{d/p}}.
\end{equation}
For the above $f_{N}\in V^{N}_{K,p}$, we know from Lemma \ref{lem 3.11} that
\begin{equation}\label{eq13}
 nc_{\rho}\bigg(\|f_N\|_{L_\nu^{p}(C_{R})}^{p}-\gamma\|f_N\|_{L_\nu^{p}(\mathbb{R}^d)}^{p}\bigg)\leq\sum\limits_{j=1}^{n}|f_N(x_{j})\nu(x_j)|^{p}
\leq  n\big(c_{\rho}\gamma+C_\rho\big)\|f_N\|_{L_\nu^{p}(\mathbb{R}^d)}^{p}
\end{equation}
holds with probability at least
\begin{equation}\label{prob}
1-A\exp\bigg(-B\frac{\gamma^{2}n c_{\rho}^{2}}{12C_\rho+2\gamma c_{\rho}}\bigg).
\end{equation}
Then, it follows from \eqref{eq11}, \eqref{eq12} and \eqref{eq13} that
\begin{equation*}
nc_{\rho}\bigg(\|f\|_{L_\nu^{p}(C_{R})}^{p}-p(1+\varepsilon)^{p-1}\varepsilon-\gamma\|f_{N}\|_{L_\nu^{p}(\mathbb{R}^d)}^{p}\bigg)
-np\Big(C^\ast+\frac{\varepsilon}{(2M)^{d/p}}\Big)^{p-1}\frac{\varepsilon}{(2M)^{d/p}}
\end{equation*}
\begin{equation}\label{h6}
\leq\sum\limits_{j=1}^{n}|f(x_{j})\nu(x_j)|^{p}\leq n\big(c_{\rho}\gamma+C_\rho\big)\|f_N\|_{L_\nu^{p}(\mathbb{R}^d)}^{p}+np\Big(C^\ast+\frac{\varepsilon}{(2M)^{d/p}}\Big)^{p-1}\frac{\varepsilon}{(2M)^{d/p}}
\end{equation} holds with the same probability as \eqref{prob}.
Since $f\in V_{K,p}(R,\delta)$, we have
\begin{equation}\label{h7}
(1-\delta)\|f\|_{L_\nu^{p}(\mathbb{R}^d)}^{p}\leq\|f\|_{L_\nu^{p}(C_{R})}^{p}.
\end{equation}
Moreover, we know  from \eqref{hc1} and Lemma \ref{lm2.5} that
\begin{equation}\label{h8}
\|f_{N}\|_{L_\nu^{p}(\mathbb{R}^d)}\leq C_{K}\big\|\big(\langle f,\widetilde{\phi}_\lambda\rangle\big)_{\lambda\in \Lambda}\big\|_{\ell_\nu^{p}}\leq B_pC_K\|f\|_{L_\nu^{p}(\mathbb{R}^d)}.
\end{equation}
Note that $\|f\|_{L_\nu^{p}(\mathbb{R}^d)}=1$. Then the sampling inequality \eqref{24} follows from \eqref{h6}-\eqref{h8}.\\

\noindent \textbf{Acknowledgement }\  
The project is partially
supported by  the  Guangxi Natural Science Foundation
(Nos. 2019GXNSFFA245012, 2020GXNSFAA159076), Guangxi Key
Laboratory of Cryptography and Information Security (No. GCIS201925), Innovation Project of School of Mathematics and Computational Science, GUET Graduate Education (No. 2022YJSCX01), Guangxi
Colleges and Universities Key Laboratory of Data Analysis and Computation.

 \end{document}